\documentclass[11pt,reqno]{amsart}
\usepackage{amssymb,amsmath,amsthm}
\usepackage{fullpage}
\usepackage{xcolor}
\usepackage{hyperref}
\usepackage{enumerate}
\usepackage{tikz}

\theoremstyle{plain}
  \newtheorem{theorem}{Theorem}
  
\theoremstyle{definition}
  \newtheorem{notation}[theorem]{Notation}
\theoremstyle{remark}
  \newtheorem{remark}[theorem]{Remark}

\renewcommand{\AA}{\mathbb{A}}
\newcommand{\GG}{\mathbb G}
\newcommand{\PP}{\mathbb P}
\newcommand{\QQ}{\mathbb Q}
\newcommand{\ZZ}{\mathbb Z}

\newcommand{\X}{\mathcal{X}}
\newcommand{\Y}{\mathcal{Y}}
\newcommand{\Z}{\mathcal{Z}}

\renewcommand{\setminus}{\smallsetminus}
\newcommand{\smat}[1]{\left(\begin{smallmatrix}#1\end{smallmatrix}\right)}
\let\hom\relax
  \DeclareMathOperator{\hom}{Hom}
\newcommand{\gp}{\mathrm{gp}}

\begin{document}
\title[No degree map]{There is no degree map for 0-cycles on Artin stacks}

\author{Dan Edidin}
\address{Department of Mathematics, University of Missouri, Columbia, MO 65211}
\email{edidind@missouri.edu}
\author{Anton Geraschenko}
\address{Department of Mathematics, Mathematics 253-37, California Institute of Technology, Pasadena, CA 91125}
\email{geraschenko@gmail.com}
\author{Matthew Satriano}
\address{Department of Mathematics, University of Michigan, Ann Arbor, MI 2074 East Hall, Ann Arbor, MI 48109--1043}
\email{satriano@umich.edu}
 \thanks{The third author was supported by an NSF postdoctoral fellowship (DMS-1103788).}

\date{\today}

\begin{abstract}
  We show that there is no way to define degrees of 0-cycles on Artin stacks with proper good moduli spaces so that (\ref{item:class-of-point}) the degree of an ordinary point is non-zero, and (\ref{item:closed-immersions}) degrees are compatible with closed immersions.
\end{abstract}
\maketitle

The theory of algebraic cycles on stacks is well-developed \cite{Gil:84,Vis:89, EdGr:98, Kre:99}. On a proper Deligne-Mumford stack, the degree of a 0-cycle is a well-defined invariant, and is compatible with proper pushforward. For quotient Deligne-Mumford stacks there is a Hirzebruch-Riemann-Roch theorem which computes the Euler characteristic of a coherent sheaf as an integral (i.e.~degree of a 0-cycle) on the associated inertia stack \cite{Toe:99, Ed:12}.

Artin stacks which are not Deligne-Mumford typically have non-proper stabilizer groups, so they are not separated and therefore not proper. However, Artin stacks with proper good moduli spaces (in the sense of \cite{Alp:08}) in many ways play the role of ``proper Artin stacks.'' For example, such stacks are universally closed \cite[Theorem 4.16(ii)]{Alp:08}, weakly separated \cite[Proposition 2.17]{ASW:10}, and satisfy formal GAGA \cite{GZB:12}. Coherent sheaves on such stacks have a finite number of non-zero cohomology groups, all of which are finite-dimensional (an easy consequence of \cite[Theorem 4.16(x)]{Alp:08}), and therefore have Euler characteristics.

It is natural to ask if there is a theory of degrees of 0-cycles on such stacks, and if there is a Hirzebruch-Riemann-Roch theorem which computes the Euler characteristic of a coherent sheaf as the degree of a 0-cycle. Another motivation for a theory of degrees is to construct numerical invariants of a variety $X$ by finding natural 0-cycles on moduli stacks associated to $X$, as is done in Gromov-Witten and Donaldson-Thomas theory. At the very least, one would like such a theory of degrees to be compatible with pushforward along closed immersions (if not proper morphisms in general), and to extend the usual notion of degree when the stack happens to be a scheme. The purpose of this note is to show there can be no such theory.

\begin{theorem}\label{thm:main}
  Let $G$ be a group, e.g.~$G=\QQ$. There do not exist homomorphisms
  \[
    \{\deg_\Y\colon A_0(\Y)\to G \mid \Y\text{ has proper good moduli space}\}
  \]
  so that
  \begin{enumerate}[{\rm (i)}]
    \item\label{item:class-of-point} the degree of a representable (i.e.~non-stacky) closed point is non-zero, and
    \item\label{item:closed-immersions} for any closed immersion $i\colon \Z\to \Y$ and any $\alpha\in A_0(\Z)$, $\deg_\Z(\alpha)=\deg_\Y(i_*\alpha)$.
  \end{enumerate}
  Specifically, over a field $k$, there do not exist such degree maps for the closed substacks of $\X = [\AA^5\setminus V(xyz,zw,v)/\GG_m^3]$, where $\GG_m^3$ acts by $(s,t,u)\cdot (x,y,z,w,v)=(sx,ty,uz,stw,stuv)$.
\end{theorem}

%

\begin{notation}\label{not:weights}
  Suppose $X\subseteq \AA^n$ is a $\GG_m^n$-invariant subscheme. Suppose a group $H$ acts on $\AA^n$ with weights $h_1,\dots, h_n\in \hom_\gp(H,\GG_m)$ (i.e.~$H$ acts via the homomorphism $(h_1,\dots, h_n)\colon H\to \GG_m^n$). We denote the stack quotient $[X/H]$ by $[X/_{\smat{h_1& \cdots& h_n}}H]$.
\end{notation}

For example, $\hom_\gp(\GG_m^3,\GG_m)\cong\ZZ^3$ and $\X = [\AA^5\setminus V(xyz,zw,v)/_{\smat{1&0&0&1&1\\ 0&1&0&1&1\\ 0&0&1&0&1}}\GG_m^3]$.

\begin{proof}
  We first show that all closed substacks of $\X$ have proper good moduli spaces. By \cite[Lemma 4.14]{Alp:08}, if $\Z\subseteq \Y$ is a closed substack and $\Y$ has good moduli space $Y$, then $\Z$ has a good moduli space which is closed in $Y$. It therefore suffices to show $\X$ has proper good moduli space.

  Consider the decomposition $\GG_m^3 = \GG_m^2\times \GG_m = \{(s,t,s^{-1}t^{-1})\}\times \{(1,1,u)\}$. We have that the ring of invariants $k[x,y,z,w,v]^{\GG_m^2}$ is $k[xyz,zw,v]$, so the induced map $[\AA^5/\GG_m^2]\to \AA^3$ is a good moduli space, on which the remaining $\GG_m$ acts with weights $\begin{pmatrix}1&1&1\end{pmatrix}$. The property of being a good moduli space morphism is local on the base in the smooth topology \cite[Proposition 4.7]{Alp:08}, so $[(\AA^5\setminus V(xyz,zw,v))/\GG_m^2]\to \AA^3\setminus 0$ is a good moduli space morphism, and hence so is
  \[
   \X = [(\AA^5\setminus V(xyz,zw,v))/\GG_m^3]\to [(\AA^3\setminus 0)/_{\smat{1&1&1}}\GG_m]=\PP^2.
  \]

  Next, note the Chow ring $A_*(\X)$ is generated by the classes of the coordinate hyperplane divisors:
  \[\begin{tabular}{c|l}
    divisor & class in $A_1(\X)$\\ \hline
    $V(x)$ & $s$\\
    $V(y)$ & $t$\\
    $V(z)$ & $u$\\
    $V(w)$ & $s+t$\\
    $V(v)$ & $s+t+u$
  \end{tabular}
  \]
  Specifically, $A_0(\X)$ is generated as an abelian group by the monomials of degree 2 in $s$, $t$, and $u$. Suppose we have a homomorphism $\deg=\deg_\X\colon A_0(\X)\to G$
  such that the induced homomorphisms $A_0(\Z)\to G$ for all closed substacks $\Z\subseteq\X$ satisfy condition (\ref{item:closed-immersions}). Replacing $G$ by the image of $\deg$, we may assume $G$ is abelian, so we write the group operation additively.

  Consider the automorphism $\alpha$ of $\X$ defined by $(x,y,z,w,v)\mapsto(y,x,z,w,v)$ on $\AA^5$ and $(s,t,u)\mapsto(t,s,u)$ on $\GG_m^3$. Note that $\alpha$ swaps the 0-cycles $V(x,z)$ and $V(y,z)$. Since automorphisms are closed immersions, condition (\ref{item:closed-immersions}) implies
  \begin{equation}\label{su=tu}
    \deg(su)=\deg(tu).
  \end{equation}

  The closed substack $V(x)$ is $[\AA^4\setminus V(zw,v)/_{\smat{0&0&1&1\\ 1&0&1&1\\ 0&1&0&1}}\GG_m^3]$. We have an automorphism $\beta$ of $V(x)$ defined by $(y,z,w,v)\mapsto(y,w,z,v)$ on $\AA^4$ and $(s,t,u)\mapsto(ut^{-1},t,st)$ on $\GG_m^3$. Since $\beta$ swaps the 0-cycles $V(x,w)$ and $V(x,z)$, we have
  \begin{equation}\label{s(s+t)=su}
    \deg (s(s+t))=\deg (su).
  \end{equation}

  Similarly, the closed substack $V(w)$ is $[\AA^4\setminus V(xyz,v)/_{\smat{1&0&0&1\\ 0&1&0&1\\ 0&0&1&1}}\GG_m^3]$. There is an automorphism $\gamma$ of $V(w)$ given by $(x,y,z,v)\mapsto(z,x,y,v)$ on $\AA^4$ and $(s,t,u)\mapsto(u,s,t)$ on $\GG_m^3$. Since $\gamma$ cyclicly permutes the 0-cycles $V(x,w)$, $V(y,w)$, and $V(z,w)$, we have
  \begin{equation}\label{s(s+t)=u(s+t)}
    \deg(s(s+t))=\deg(t(s+t))=\deg(u(s+t)).
  \end{equation}

  By (\ref{s(s+t)=su}) and (\ref{s(s+t)=u(s+t)}), we see that $\deg(su+tu)=\deg(su)$, so $\deg(tu)=0$. Then by (\ref{su=tu}), $\deg(su)=0$. Combining this with (\ref{s(s+t)=su}) and (\ref{s(s+t)=u(s+t)}), we have
  \[
    \deg((s+t)(s+t+u)) \overset{\mbox{\scriptsize(\ref{s(s+t)=u(s+t)})}}{=} 3\deg((s+t)s) \overset{\mbox{\scriptsize(\ref{s(s+t)=su})}}{=} 3\deg(su) = 0.
  \]

  The closed substack $V(w,v)$ of $\X$ is $[\AA^3\setminus V(xyz)/_{\smat{1&0&0\\ 0&1&0\\ 0&0&1}}\GG_m^3]$, a non-stacky closed point whose class in $A_0(\X)$ is $(s+t)(s+t+u)$, so the above equation contradicts condition (\ref{item:class-of-point}).
\end{proof}

\begin{remark}
  We believe condition (\ref{item:closed-immersions}) alone is not enough to determine the degree of every element of $A_0(\X)$. Specifically, it seems the degree of the 0-cycle $V(x,y)\subseteq \X$, which has class $st\in A_0(\X)$, is not determined. Geometrically, this 0-cycle is a curve with 1-dimensional generic stabilizer.
\end{remark}

\begin{remark}
  The stack $\X$ is a toric stack (in fact, a fantastack) in the sense of \cite[\S4]{toricartin1}. Following \cite[Notation 4.8]{toricartin1}, the stacks $\X$, $V(x)$, and $V(w)$ are depicted below. This point of view makes it easy to see the existence of the automorphisms $\alpha$, $\beta$, and $\gamma$.
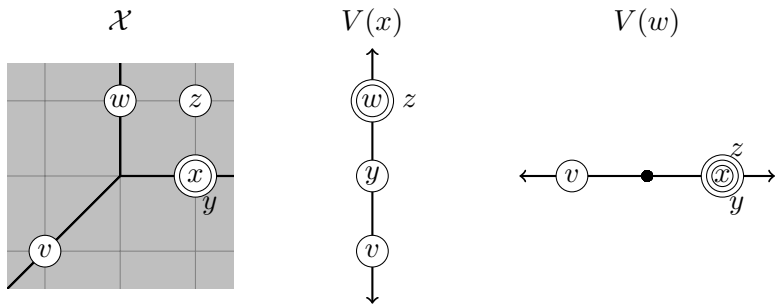
  \[\begin{tabular}{ccc}
   $\X$ & $V(x)$ & $V(w)$ \\
   \begin{tikzpicture}
   \clip (-1.5,-1.5) rectangle (1.5,1.5);
   \foreach \x/\y/\z/\w in {9/0/0/9, 0/9/-9/-9, -9/-9/9/0}
     \filldraw [draw=black,fill=lightgray] (\x,\y) -- (0,0) -- (\z,\w);
   \draw[help lines] (-10,-10) grid (10,10);
   \foreach \x/\y/\z/\w in {9/0/0/9, 0/9/-9/-9, -9/-9/9/0}
     \draw [thick] (\x,\y) -- (0,0) -- (\z,\w);
   \filldraw[draw=black,fill=white] (1,0) node[label=below:$\quad y$] {} circle (8pt);
   \foreach \x/\y/\dottext in {1/0/x, 1/1/z, 0/1/w, -1/-1/v}
     \filldraw[draw=black,fill=white] (\x,\y) node {$\dottext$} circle (6pt);
   \end{tikzpicture}
   &
   \raisebox{-3mm}{\begin{tikzpicture}
     \clip (-1.5,-1.8) rectangle (1.5,1.8);
     \draw[thick,<->] (0,-1.7) -- (0,1.7);
     \filldraw[draw=black,fill=white] (0,1) node[label=right:$\ z$] {} circle (8pt);
     \filldraw[draw=black,fill=white] (0,1) node {$w$} circle (6pt);
     \filldraw[draw=black,fill=white] (0,0) node {$y$} circle (6pt);
     \filldraw[draw=black,fill=white] (0,-1) node {$v$} circle (6pt);
   \end{tikzpicture}}
   &
   \begin{tikzpicture}
     \clip (-1.8,-1.5) rectangle (1.8,1.5);
     \draw[thick,<->] (-1.7,0) -- (1.7,0);
     \draw plot[only marks,mark=*] (0,0);
     \filldraw[draw=black,fill=white] (1,0) node[label=below:$\quad y$,label=above:$\quad z$] {} circle (8pt);
     \filldraw[draw=black,fill=white] (1,0) node {} circle (6pt);
     \filldraw[draw=black,fill=white] (1,0) node {$x$} circle (4pt);
     \filldraw[draw=black,fill=white] (-1,0) node {$v$} circle (6pt);
   \end{tikzpicture}
   \end{tabular}
  \]
\end{remark}

\textbf{Acknowledgements.} We would like to thank Tom Graber, who got the second and third authors interested in this problem, and skillfully disposed of many proposed definitions of degree.


\end{document}